# A Constructive Proof of the Four Color Theorem

**Author: Ding Dagong**

**Abstract** This paper presents an approach to proving the Four Color Theorem that differs from traditional "reducible configuration" methods. By introducing concepts such as "simple islands", "outer boundaries", "primitive sets", "Property A", "knots", and "valid number-pair groups", together with the operation of "adding an n-point region to an interval", we show that if every simple island can be colored with four colors, then any given map can be colored with four colors. We construct a framework for coloring an arbitrarily given spherical map step by step. At the heart of coloring a simple island are four theorems, which ensure that after successively adding specific regions to an outer boundary satisfying Property A, the new outer boundary still satisfies Property A, so that the simple island can ultimately be colored with four colors, and in turn the entire arbitrarily given spherical map can be colored. This method avoids computer enumeration and offers a more constructive perspective on the proof.

## Introduction

The Four Color Theorem asserts that any spherical map can be colored with four colors so that adjacent regions receive different colors. Since it was proposed in 1852, the proof of this theorem has undergone many twists and turns. Although Appel and Haken proved the theorem with the aid of computers in 1976, their traditional proof relies on the enumeration of a large number of reducible configurations and lacks intuitiveness and constructiveness. Inspired by the idea of Kempe chains but dispensing with complex configuration analysis, this paper proposes a straightforward constructive method: for any given spherical map, color it with four colors through a series of explicit operational steps. To this end, we introduce concepts such as "island", "outer boundary" (outer boundary), "primitive set" (primitive set), "Property A", and "knot", and establish three key theorems to guarantee the feasibility of the coloring process.

**Main Text**

Given an arbitrary spherical map, we assume that every node in the given map has degree 3. If the given map has several nodes of degree greater than 3, we may add a small region over each such node to cover it, so that all nodes of the map become degree 3. If the map with the small regions added can be colored with four colors, then after coloring we remove these added small regions; the colors of the regions in the original map remain unchanged, because regions that were adjacent in the original map remain adjacent after the small regions are added, and therefore all regions of the original map can be colored with four colors.

We now set out to color an arbitrarily given spherical map. Since the map is given, its topological structure is completely determined and unique.

We first choose an arbitrary region in the given spherical map and color it. After it is colored, among the three regions meeting at each node of the map, either one is already colored or none of the three is colored. At a node where one region is already colored, two of the edges form part of the boundary while the other does not; hence every node on the boundary has two edges belonging to the boundary, and consequently the boundary in the given map can only consist of several simple closed curves. We call the chosen colored region the sea, and each simple closed curve together with all the uncolored regions inside it an island belonging to this sea. This sea is also called the sea belonging to that island.

Choose an arbitrary point on the sea and call it the pole. If on the same island there are two distinct regions A and B, and every connection from the pole to region B must pass through region A, then we say that region B lies inside region A, or that region B is an interior region of region A and region A is an exterior region of region B.

We now remove all interior regions from one of the islands in the given map. This yields an island with no interior regions, because removing an interior region does not affect the adjacency relations among its exterior regions and other regions, so the regions that remain keep the same adjacency relations as before. An island with no interior regions is called a simple island.

We now show that if, for a colored sea and any simple island belonging to that sea, we can color all regions of the simple island with four colors, then we can color any given spherical map.

We first choose an arbitrary island in the sea and remove its interior regions, obtaining a simple island. By hypothesis, every region on the island can be colored with one of the four colors. Now we regard each region of this simple island as a distinct, already colored sea. We then put back the interior regions previously removed from the island at their original positions; these restored, uncolored regions become regions on islands of these new seas.

We repeat this procedure: choose an arbitrary island and the sea belonging to it among the old and new seas, remove all interior regions of the island, then color the simple island stripped of its interior regions, thereby obtaining further new seas. We continue until none of the old or new seas contains any island. This is achievable because we are dealing with a finite map, and each application of the procedure always adds some colored regions as new seas and reduces the number of uncolored regions, so the process must terminate. That is, all regions in the arbitrarily given map are colored.

We now turn to the problem of coloring the regions of a simple island in the sea.

## Preliminaries

### 1.1 Basic Setup

A simple closed curve in a spherical map is called an outer boundary; it divides the sphere into two planes. During coloring, one of the planes has all of its regions colored; we call this plane the inside of the boundary, while the other plane, whose regions are all uncolored, is called the outside. Because every node has degree 3, each node on the outer boundary has one further edge, pointing either into the inside or into the outside. We call a node whose additional edge points outside an outer point of the outer boundary, and a node whose additional edge points inside an inner point of the outer boundary. A region with no inner points is simply a polygon. Two adjacent outer points on the outer boundary, together with the arc of the outer boundary between them, are called an interval. Coloring an interval, as discussed below, means coloring the region outside the boundary that the interval represents.

If we are faced with a colored sea and an island in that sea, then because the sea is already colored, the sea is the interior region of the island's coastline, which is a simple closed curve. On this closed curve there are no inner points, only outer points. Hence the sea together with the island's coastline can be regarded as a colored polygon.

### 1.2 Primitive Set

Given an outer boundary, consider all possible colorings of the regions inside the boundary that satisfy the four-color condition. For each such coloring, take the set of all colorings of the intervals on the outer boundary that satisfy the four-color condition; the union of all these sets is called the primitive set of that outer boundary. The primitive set of an outer boundary is the basis for all subsequent coloring operations on regions outside the boundary.

### 1.3 Region-Adding Operation

Adding an n-point region to an interval of an outer boundary means adding a region that connects to the outer boundary only through that interval and has n nodes not lying on the original outer boundary. Once the new region is added, a new outer boundary is formed, and the added region becomes an interior region of the new outer boundary, its color being the color of that interval. The n nodes not on the original outer boundary are called new points, and these n new points are also outer points of the new outer boundary. Points that were outer points of the original outer boundary remain outer points of the new outer boundary. On the new outer boundary, an interval whose two endpoints are both old points is called an old interval; an interval whose two endpoints are both new points is called a new interval; and an interval with one old and one new endpoint is called a junction interval. After adding an n-point region ($n>0$), the new outer boundary has two junction intervals.

Adding a 0-point region and adding a 1-point region are jointly referred to as adding a simple region. Adding a simple region creates no new intervals, and it reduces the number of intervals on the outer boundary. Adding a 0-point region to an interval merges the intervals on the two sides of that interval into a single interval in the new outer boundary, so those two intervals must have the same color. Adding a 1-point region to an interval makes the intervals on its two sides adjacent in the new outer boundary, so those two intervals must have different colors. In fact, if there is a chain of consecutive intervals with more than two intervals, one can always add some simple regions within the chain so that the two end intervals of the chain either become a single interval or become adjacent intervals.

When the outer boundary has three intervals, their colors must be pairwise different, so a 1-point region can be added to any one of them. When the outer boundary has two intervals, the intervals on the two sides of each interval are actually the same interval, so a 0-point region can be added to either interval.

### 1.4 Property A

Suppose the primitive set of an outer boundary is nonempty. We choose an arbitrary interval on this outer boundary, add a simple region to it, then add a simple region to an arbitrary interval on the resulting outer boundary, and repeat this operation until only two or three intervals remain on the outer boundary. Note that each addition of a simple region imposes a requirement on the coloring of the intervals of the previous outer boundary. If in the primitive set of the initial outer boundary there always exists an element satisfying the requirement of each such arbitrary addition of a simple region, we say that the primitive set of this outer boundary satisfies Property A, or simply that the outer boundary satisfies Property A.

The Four Color Theorem requires that adding an arbitrary n-point region to an interval of an outer boundary must still allow coloring with four colors, whereas satisfying Property A only requires that adding simple regions to the intervals of the outer boundary allows coloring with four colors. Hence an outer boundary satisfying Property A is a necessary condition for the Four Color Theorem to hold.

## 2. Main Theorems

**Theorem 1** When coloring an island of the sea, the sea, regarded as a region, has an outer boundary that satisfies Property A. (Fundamental Theorem of Polygons)

**Proof** We are faced with a spherical map consisting of a colored sea and an uncolored island in that sea. In this spherical map, the sea is simply a colored polygon.

For a digon or a triangle, the theorem clearly holds. Assume that for $m < n$, the outer boundary of an m-gon satisfies Property A.

For an n-gon, if a 1-point region is added to an interval of its outer boundary, a triangle appears in the map. We know that when coloring a map with four colors, a triangle can be reduced to a single point for coloring, because if the reduced map can be colored with three colors, then restoring the point to a triangle, the triangle can always be colored with the fourth color, different from the three colors of the three adjacent regions.

If a 0-point region is added to its outer boundary, a digon appears in the map. By reasoning similar to the above, the digon together with its two incident edges can be reduced to a single edge for coloring.

Therefore the coloring of the intervals on the outer boundary obtained after adding a simple region to a polygon corresponds to the coloring of the intervals on the outer boundary of an (n-2)-gon or an (n-1)-gon. By the induction hypothesis, the outer boundary of the n-gon satisfies Property A. Hence the outer boundary of any polygon satisfies Property A, that is, the outer boundary of the sea satisfies Property A.

**Theorem 2** If an outer boundary satisfies Property A, then after adding a 2-point region to any interval of that outer boundary, the resulting new outer boundary also satisfies Property A. (2-Point Region Theorem)

**Proof** Figure 1 shows the map obtained by adding a 2-point region to the interval BC of an outer boundary; adding the 2-point region creates only one new interval MN. If 1-point regions are repeatedly added to the junction intervals of Figure 2, Figure 4 is obtained. We will explain the existence of Figure 4 later. The interval ST in Figure 4 is in fact an extension of the interval MN in Figure 1; we also call the interval ST in Figure 4 a new interval, and the points S and T new points.

If we add a 0-point region to the junction interval AM in Figure 1 (see Figure 2), or to the junction interval FS in Figure 4 (see Figure 3), Figures 2 and 3 can also be regarded as obtained by first adding a 1-point region to the interval AB or FG, respectively, on the original outer boundary corresponding to the junction interval AM in Figure 1 or FS in Figure 4, and then adding a 1-point region, or repeatedly adding 1-point regions, to the junction interval. Because the original outer boundary satisfies Property A, both Figures 2 and 4 satisfy Property A.

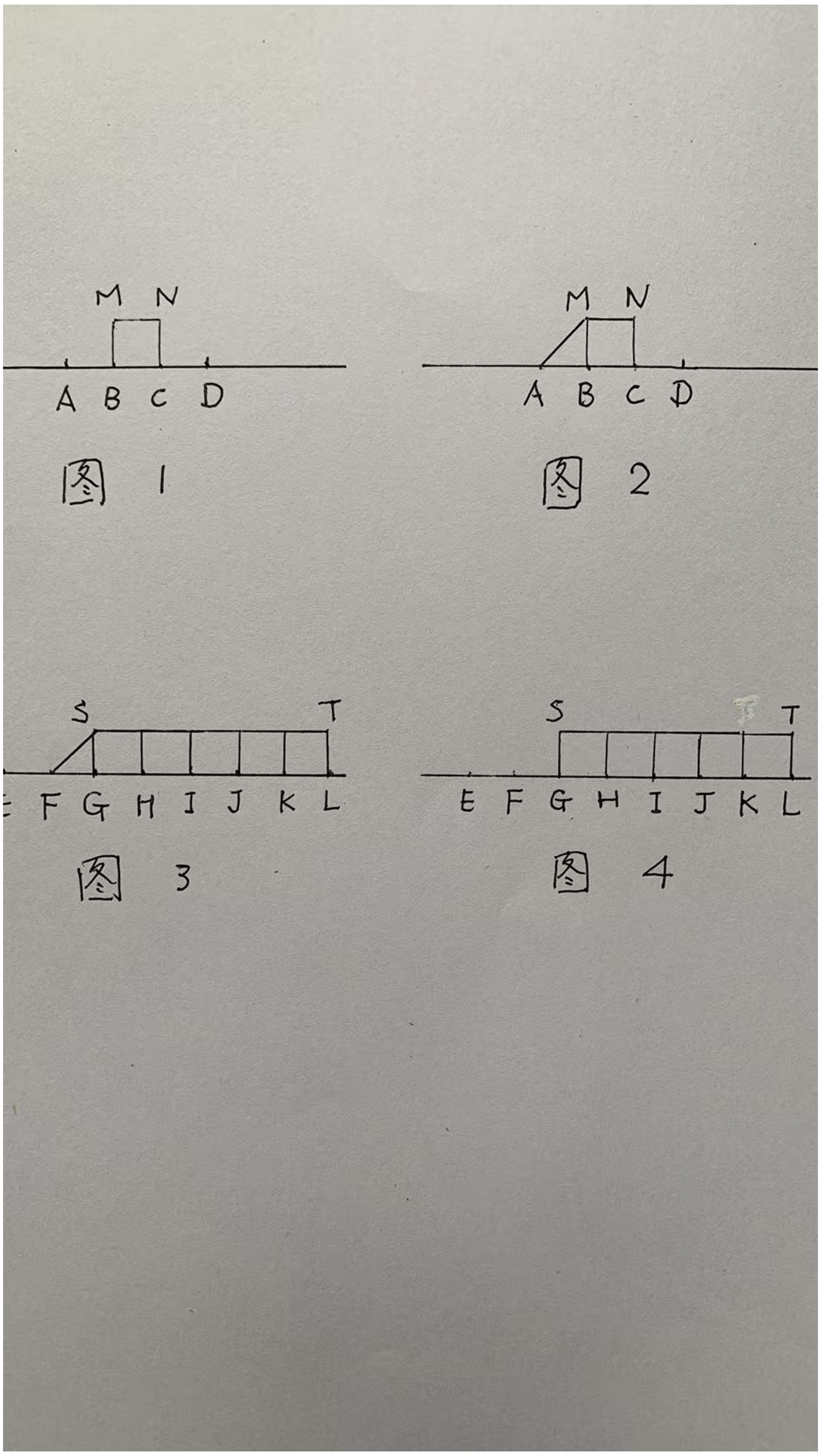

If a simple region is added to an old interval in Figure 1 or Figure 4, since the old interval lies on the original outer boundary satisfying Property A, this simple region can always be added, and the problem remains the old one, namely, "does adding a 2-point region to any interval satisfying Property A still yield Property A?" Only the number of old intervals has decreased.

Regarding the existence of Figures 1 and 4, we now consider the problem from another angle. If a node is incident to three edges, it is of course incident to three regions, and these three regions are pairwise adjacent, so their colors are pairwise different. We label the four colors 1, 2, 3, and 4. Imagine a tetrahedron whose four faces are painted with four different colors (see Figure 5). The three colors of the regions meeting at each node in the colored map correspond to three mutually distinct numbers. On the tetrahedron in Figure 5, find the vertex where these three numbers meet; if, viewed from outside, the three numbers of the regions meeting at this node in the map are arranged in the same direction as the three numbers at that vertex in Figure 5, then the value of this node is defined to be 1. For example, suppose the colors of the three regions meeting at a certain node are 1, 2, and 3. Viewing the tetrahedron in Figure 5 from outside, the numbers 1, 2, 3 are arranged counterclockwise at the vertex where 1, 2, 3 meet; if the colors 1, 2, 3 at that node in the map are also arranged counterclockwise, the value of the node in the map is defined to be 1; otherwise it is defined to be -1. The numbers 1, -1, and also 0 mentioned here are all numbers modulo 3. For a polygon, suppose its color is labeled 4; then the colors of the intervals on its outer boundary are 1, 2, and 3. Without loss of generality, start from the interval of color 1 and examine the colors of the intervals one after another in the same direction. This is in effect examining the values of the vertices one after another, where a vertex takes only the value 1 or -1. If the values 1 and -1 occur equally often, or one occurs a multiple of 3 more than the other, then when we return to the initial interval we return to the original first color. In that way the polygon and the polygons adjacent to it can be colored. Otherwise, on returning to the initial interval we do not recover the original first color, which is not allowed. Hence, for every polygon in the given map, the sum of the values assigned to all its vertices (which can only be 1 or -1) being 0 (mod 3) is a necessary and sufficient condition for the given map to be colorable according to the four-color rule.

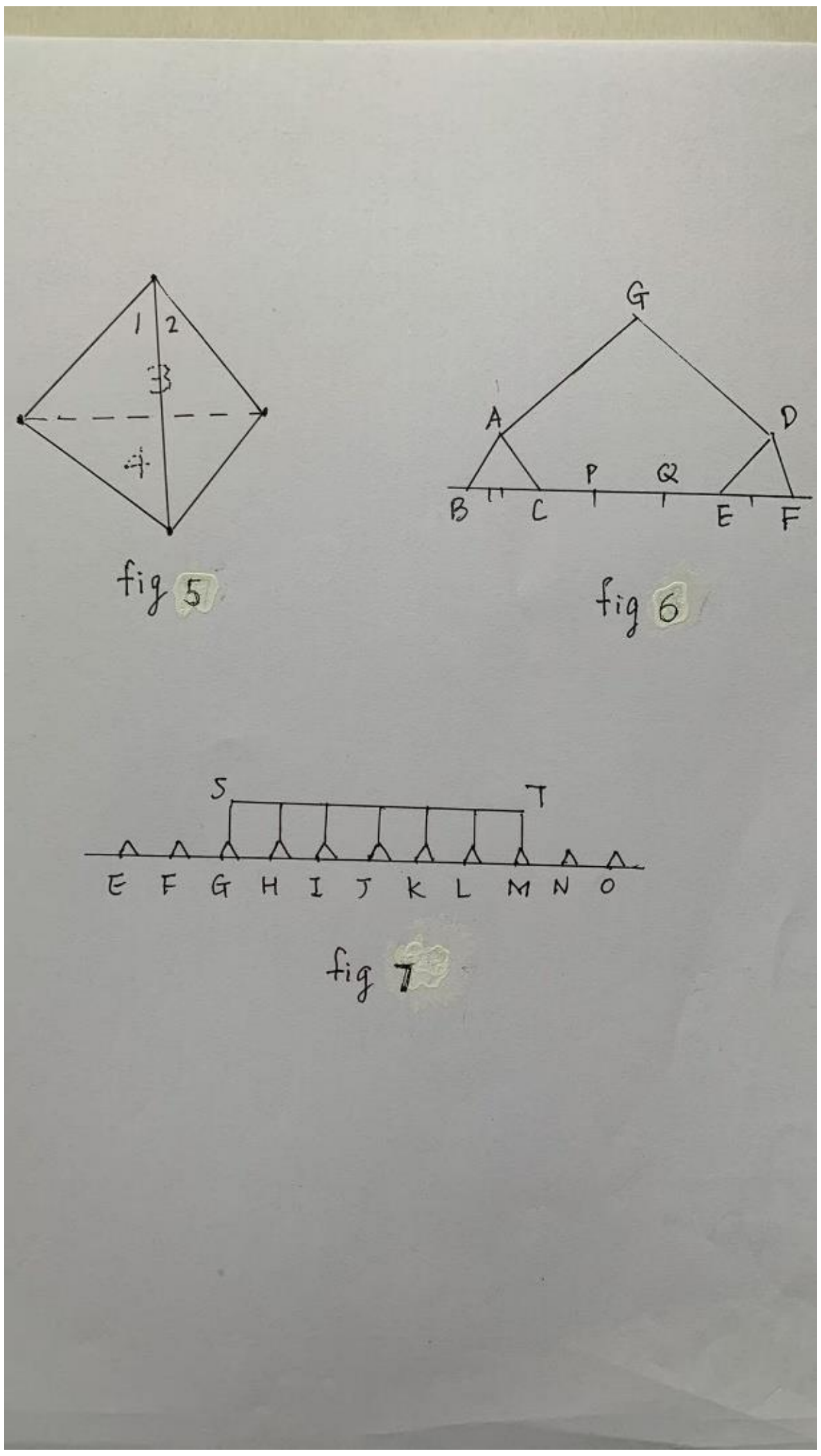

We now introduce a new concept: the "knot".

For a polygon, every node on its outer boundary is an outer point, and each of these outer points can be regarded as a knot. Now the outer boundary has only knots. Two adjacent knots, together with the arc of the outer boundary between them, are called an interval. Note: as we shall see later, under this concept two adjacent intervals may overlap over a small arc. If a 0-point region is added to an interval, the two knots at the ends of that interval cease to be knots, and the outer boundary has two fewer knots. If a 1-point region is added to an interval, we call the added new point together with that interval a new knot. The two knots formerly at the ends of the interval are now part of the new knot, so the number of knots on the outer boundary decreases by one. The newly added point is called the vertex of this new knot. Adding a 1-point region to two adjacent intervals amounts to adding a 1-point region at the vertices of these two knots. Of course, there may still be some inner points between the two knots at the ends of the interval where the 1-point region is added (see Figure 6). In Figure 6, knots ABC and DEF are two adjacent knots on the original outer boundary, knot GBF is the new knot formed after adding a 1-point region to knots ABC and DEF, and points P and Q are inner points between two adjacent knots on the original interval.

Each knot contains a number of number-pairs. A number-pair represents a contribution that can be made simultaneously to its left interval and right interval. For example, a node, regarded as a knot, can contribute the value 1, or -1, simultaneously to the interval on its left and the interval on its right. Thus a node as a knot has two number-pairs, namely (1,1) and (-1,-1). We call the following two types of number-pair groups valid number-pair groups. A valid number-pair group of the first type contains two number-pairs whose two left values are different, and whose right values are simultaneously one greater than their left values, or simultaneously one less than their left values, or simultaneously equal to their left values. For example, the two number-pairs of the node just discussed form a valid number-pair group. A valid number-pair group of the second type contains three number-pairs whose three left numbers take exactly two values, while of the three right numbers one is one greater than its left value, one is one less than its left value, and one is equal to its left value. For example: (1,0); (1,-1); (0,0). These three number-pairs satisfy the condition for a valid number-pair group. If the three left numbers took a single value, then requiring the three right numbers to be one greater, one less, and one equal would make the three right numbers take three distinct values; but here the three left numbers take exactly two values, and it is easy to verify that in this case the three right numbers also take exactly two values. We now show that adding a 1-point region between two knots each containing a valid number-pair group yields a new knot that also contains a valid number-pair group.

We choose an arbitrary number-pair from the left valid number-pair group and an arbitrary number-pair from the right valid number-pair group. The left value of the number-pair chosen from the left knot's valid number-pair group is called a left value of this choice; the right value of the number-pair chosen from the right knot's valid number-pair group is called a right value of this choice; and the sum of the right value of the number-pair chosen from the left knot's group plus the left value of the number-pair chosen from the right knot's group is called an inner value of this choice. Left value, right value, and inner value are three numbers associated with a particular choice of number-pairs. If the value of the vertex is 1, then left value plus 1, right value plus 1, and inner value plus 1 are called the plus-1 left sum, plus-1 right sum, and plus-1 inner sum, respectively. Similarly, when the value of the vertex is -1, we have the minus-1 left sum, minus-1 right sum, and minus-1 inner sum. The plus-1 inner sum and the minus-1 inner sum are collectively called inner sums. The plus-1 left sum and plus-1 right sum form a number-pair, which we call the plus-1 pair; the minus-1 left sum and minus-1 right sum form a number-pair, which we call the minus-1 pair.

Left value, right value, inner value, plus-1 left sum, plus-1 right sum, plus-1 inner sum, minus-1 left sum, minus-1 right sum, minus-1 inner sum, inner sum, plus-1 pair, and minus-1 pair are all concepts defined with respect to a particular way of choosing two number-pairs. In addition, as mentioned earlier, there may be some inner points between the two knots where a 1-point region is added; we call the sum of the values of all these inner points the inner-point sum.

Since the map can be colored with four colors, the sum of the values at all nodes of any polygon must be 0, so the inner-point sum plus the inner sum must equal 0. All minus-1 pairs when the inner value is 1 (in which case the vertex has value -1), together with all plus-1 pairs when the inner value is -1 (in which case the vertex has value 1), have inner sum 0. All plus-1 pairs when the inner value is 0, together with all minus-1 pairs when the inner value is -1, have inner sum 1. All minus-1 pairs when the inner value is 0, together with all plus-1 pairs when the inner value is 1, have inner sum -1.

If the inner-point sum is A, then the inner sum must be -A, where -A may be any of the three numbers 1, -1, or 0. We therefore need to show that among all number-pairs whose inner sum satisfies this requirement, there is contained a valid number-pair group. Let us first examine some properties of valid number-pair groups.

It is easy to verify that adding 1 or subtracting 1 from the left value of every number-pair in a valid number-pair group, or adding 1 or subtracting 1 from the right value of every number-pair in a valid number-pair group, yields a valid number-pair group of the same type.

For a valid number-pair group consisting of two number-pairs, we can add 1,

subtract 1, or leave unchanged the two left numbers simultaneously, and do the same to the two right numbers, so that the resulting two number-pairs become (0,0) and (1,1). It suffices to add a number on the left so that the two left numbers become 0 and 1, and add a number on the right so that the two right numbers also become 0 and 1.

For a valid number-pair group consisting of three number-pairs, we can add 1, subtract 1, or leave unchanged the three left numbers simultaneously, and do the same to the three right numbers, so that the resulting three number-pairs become (0,0), (0,1), and (1,0), or become (0,0), (0,-1), and (-1,0). It suffices to add a number on the left so that the two left numbers with equal values become 0, and add a number on the right so that the two right numbers with equal values also become 0. For the three number-pairs to form a valid number-pair group, the two remaining nonzero numbers must be equal.

In addition, if we swap the right and left values of every number-pair in a valid number-pair group, the result is a valid number-pair group of the same type.

If we add the same number to the right number of every number-pair in the left knot, we may regard this number as taken from the inner-point sum. In this case, to prove that the new knot contains a valid number-pair group, we still need to show that among all number-pairs whose inner sum equals any of the three numbers 0, 1, -1, there is contained a valid number-pair group. Similarly, adding the same number to the left number of every number-pair in the right knot, a number we may likewise regard as taken from the inner-point sum, requires us to show that among all number-pairs whose inner sum takes any value, there is contained a valid number-pair group.
If we add a number to the left number of every number-pair in the left knot and a number to the right number of every number-pair in the right knot, and if in the end all number-pairs satisfying the requirement contain a valid number-pair group, we can subtract from the left side of each number-pair in that group the number added on the left, and subtract from the right side the number added on the right. The resulting valid number-pair group is precisely the valid number-pair group on the new knot after the 1-point region is added.

Because we can always add a number to the left of every number-pair and a number to the right of every number-pair of a valid number-pair group so that the final result is one of (0,0), (1,1); (0,0), (0,1), (1,0); and (0,0), (0,-1), (-1,0), we need only consider these three valid number-pair groups.

For each of these three results, swapping the left and right values of its number-pairs yields a number-pair group identical to the original valid number-pair group. Moreover, in analyzing whether a new number-pair group contains a valid number-pair group, exchanging the left and right valid number-pair groups gives the same result. Hence, if we label these three valid number-pair groups 1, 2, and 3, we need only tabulate and analyze the six cases in which the left and right labels are respectively 1,1; 1,2; 1,3; 2,2; 2,3; and 3,3.

The six tables below each correspond to one of the six cases above, and each table is divided into left, middle, and right parts. The left part of a table shows the valid number-pair group contained in the left knot and the valid number-pair group contained in the right knot when a 1-point region is added.

In the middle part of a table, the central number denotes the inner value. On its left and right are respectively the left values and right values of all choices that give this inner value. The two numbers to the left of the left values and to the right of the right values are the plus-1 pairs for the same choices. The two numbers flanking the plus-1 pairs are the minus-1 pairs for the same choices. In the right part of the table, the number on the left denotes the inner sum, and to its right are all number-pairs satisfying this inner sum. When the inner sum is 0, these are all minus-1 pairs for inner value 1 and all plus-1 pairs for inner value -1. When the inner sum is 1, these are all plus-1 pairs for inner value 0 and all minus-1 pairs for inner value -1. When the inner sum is -1, these are all minus-1 pairs for inner value 0 and all plus-1 pairs for inner value 1. We find one valid number-pair group among all number-pairs for inner sums 0, 1, and -1, respectively, and mark each number-pair of that group with an asterisk (*).

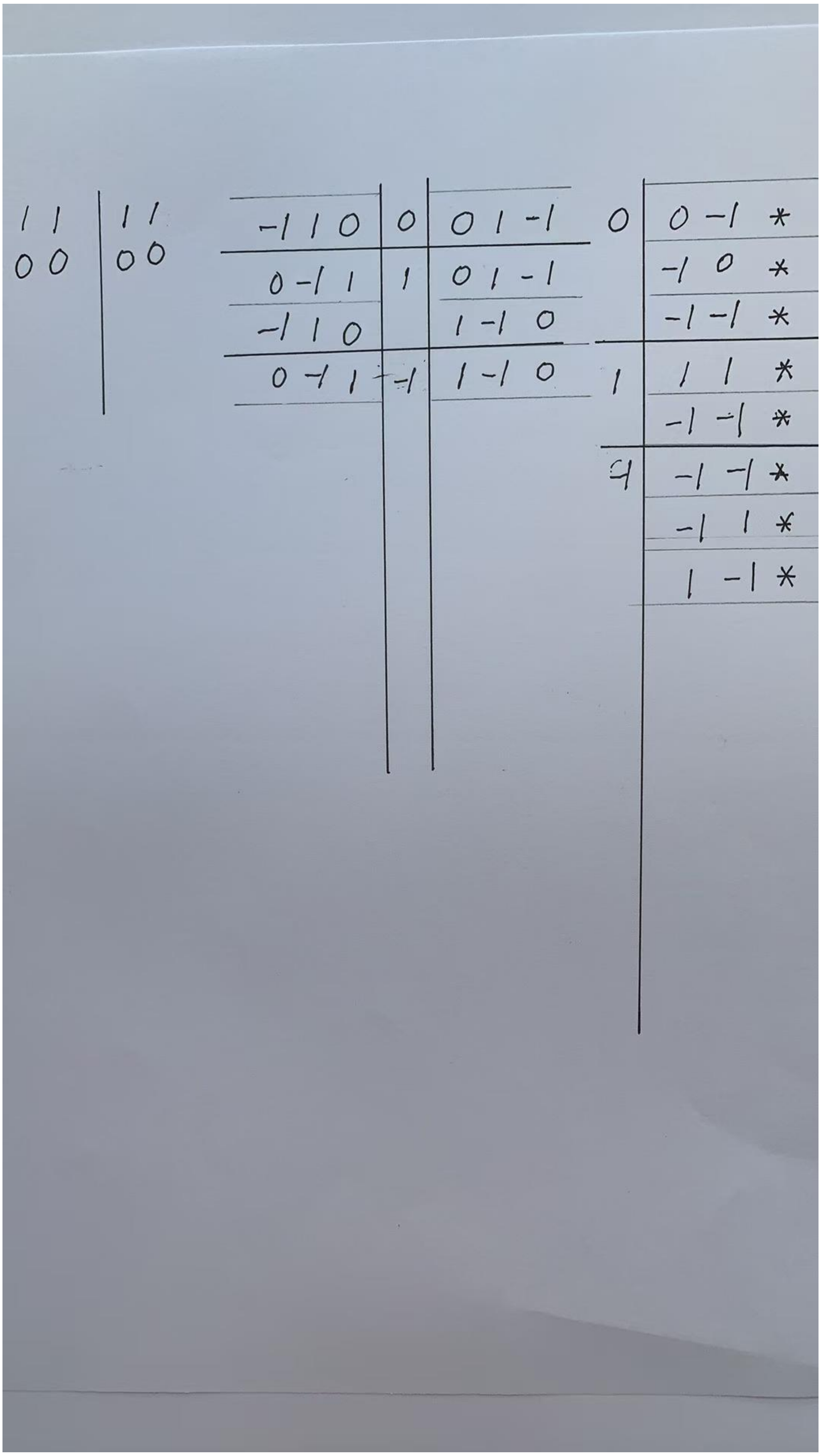

| | |
|---|---|
| 1 1 | 1 1 |
| 0 0 | 0 0 |

| | | |
|---|---|---|
| -1 1 0 | 0 | 0 1 -1 |
| 0 -1 1 | 1 | 0 1 -1 |
| -1 1 0 | | 1 -1 0 |
| 0 -1 1 | -1 | 1 -1 0 |
| 0 | 0 -1 | * |
| | -1 0 | * |
| | -1 -1 | * |
| 1 | 1 1 | * |
| | -1 -1 | * |
| -1 | -1 -1 | * |
| | -1 1 | * |
| | 1 -1 | * |

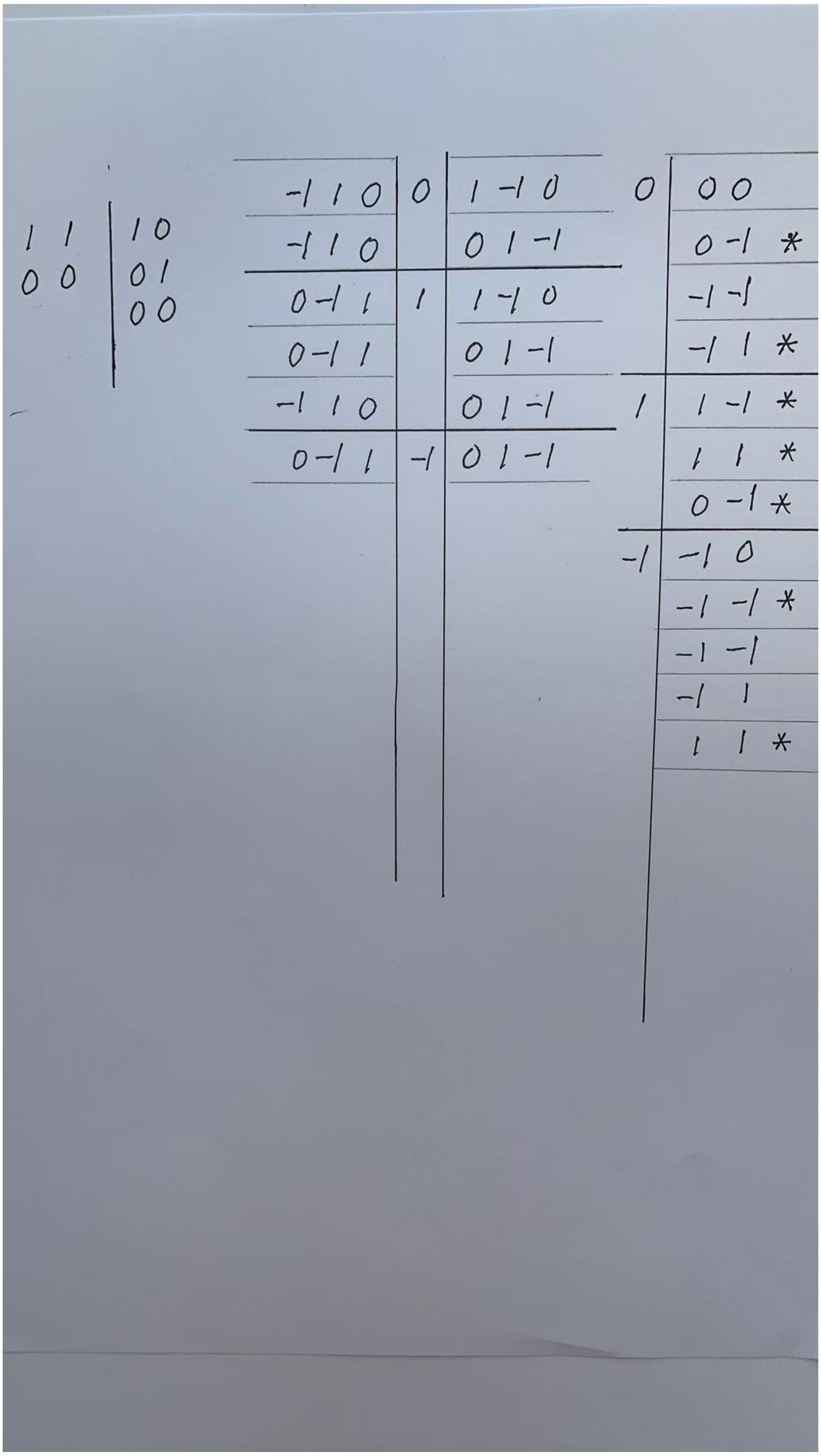

| | |
|---|---|
| 1 1 | 10 |
| 0 0 | 01 |
| | 00 |

| | | |
|---|---|---|
| -1 1 0 | 0 | 1 -1 0 |
| -1 1 0 | | 0 1 -1 |
| 0 -1 1 | 1 | 1 -1 0 |
| 0 -1 1 | | 0 1 -1 |
| -1 1 0 | | 0 1 -1 |
| 0 -1 1 | -1 | 0 1 -1 |

| | |
|---|---|
| 0 | 0 0 |
| | 0 -1 * |
| | -1 -1 |
| | -1 1 * |
| 1 | 1 -1 * |
| | 1 1 * |
| | 0 -1 * |
| -1 | -1 0 |
| | -1 -1 * |
| | -1 -1 |
| | -1 1 |
| | 1 1 * |
| 1 1 | -1 0 |
| 0 0 | 0 -1 |
| | 0 0 |

| | | |
|---|---|---|
| 0 -1 1 | 0 | 0 1 -1 |
| -1 1 0 | | -1 0 1 |
| -1 1 0 | | 0 1 -1 |
| 0 -1 1 | 1 | -1 0 1 |
| 0 -1 1 | | 0 1 -1 |
| -1 1 0 | -1 | 0 1 -1 |

| | |
|---|---|
| 0 | 0 1 * |
| | 0 -1 * |
| | 1 1 * |
| 1 | -1 1 * |
| | 1 0 * |
| | 1 1 |
| | -1 -1 |
| -1 | 0 -1 * |
| | -1 1 * |
| | -1 -1 |
| | -1 0 |
| | -1 1 |
| 1 0 | 1 0 |
| 0 1 | 0 1 |
| 0 0 | 0 0 |

| | | |
|---|---|---|
| 0 -1 1 | 0 | 1 -1 0 |
| 0 -1 1 | | 0 1 -1 |
| -1 1 0 | | 1 -1 0 |
| -1 1 0 | | 0 1 -1 |
| 0 -1 1 | 1 | 0 1 -1 |
| -1 1 0 | | 1 -1 0 |
| -1 1 0 | | 0 1 -1 |
| -1 1 0 | | 0 1 -1 |
| -1 1 0 | -1 | 0 1 -1 |

| | |
|---|---|
| 0 | 0 -1 |
| | -1 0 |
| | -1 -1 |
| | -1 -1 * |
| | 1 1 * |
| 1 | -1 -1 |
| | -1 1 |
| | 1 -1 |
| | 1 1 * |
| | -1 -1 * |
| -1 | 0 0 * |
| | 0 -1 |
| | -1 0 |
| | -1 -1 * |
| | -1 1 |
| | 1 -1 |
| | 1 1 |
| | 1 1 |
| 1 0 | -1 0 |
| 0 1 | 0 -1 |
| 0 0 | 0 0 |

| | | |
|---|---|---|
| 0 -1 1 | 0 | -1 0 1 |
| 0 -1 1 | | 0 1 -1 |
| -1 1 0 | | 0 1 -1 |
| -1 1 0 | | -1 0 1 |
| -1 1 0 | | 0 1 -1 |
| -1 1 0 | 1 | -1 0 1 |
| -1 1 0 | | 0 1 -1 |
| 0 -1 1 | -1 | 0 1 -1 |
| -1 1 0 | | 0 1 -1 |

| | |
|---|---|
| 0 | -1 1 |
| | -1 -1 * |
| | -1 1 |
| | 1 1 * |
| 1 | -1 0 * |
| | -1 1 * |
| | 1 1 * |
| | 1 0 |
| | 1 1 |
| | 0 -1 |
| | -1 -1 |
| -1 | 0 1 * |
| | 0 -1 * |
| | -1 -1 * |
| | -1 1 |
| | -1 -1 |
| | 1 0 |
| | 1 1 |

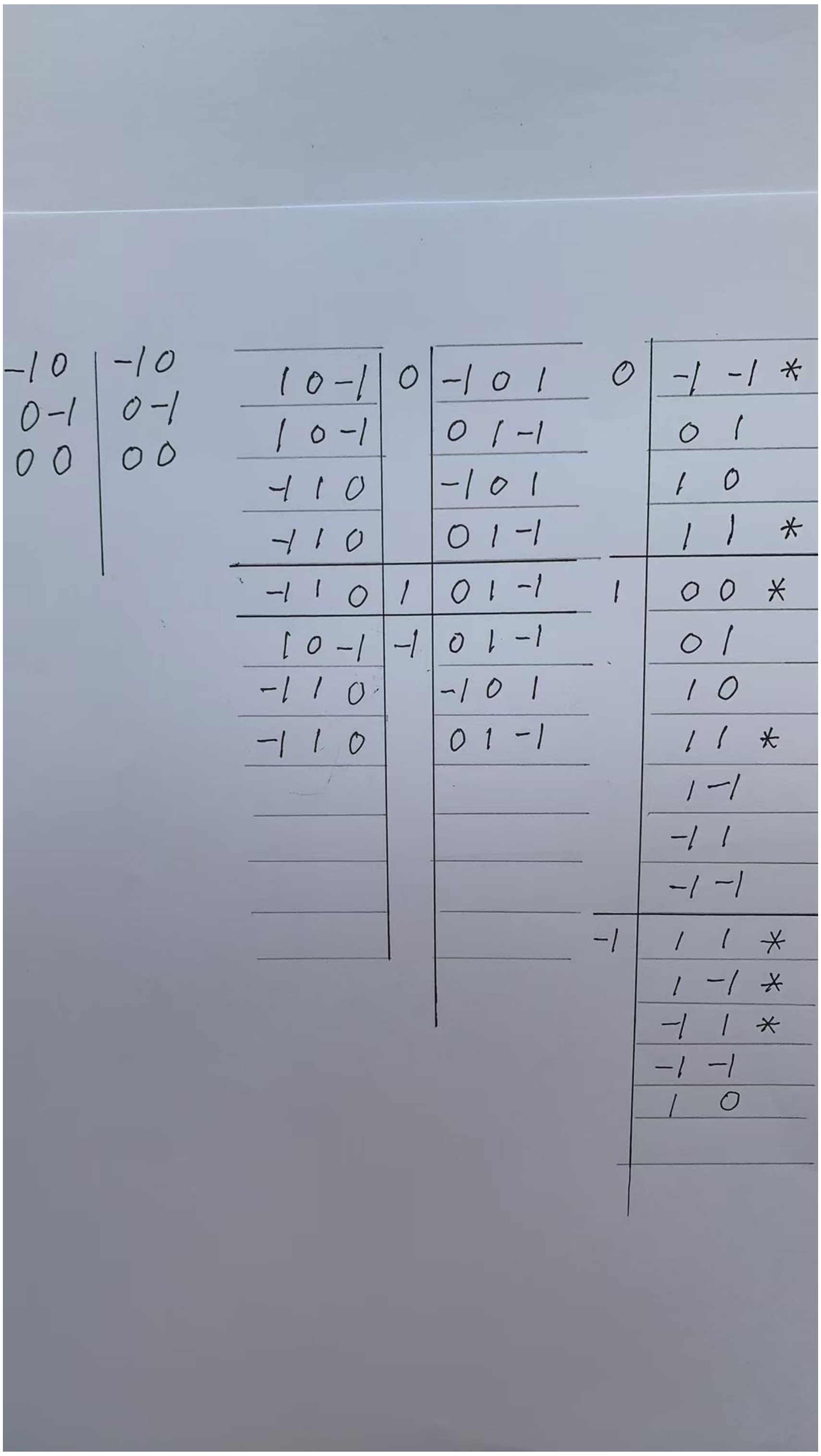

From these tables it can be seen that, whatever the inner-point sum and whatever valid number-pair groups are contained in the two knots on either side of the interval, adding a 1-point region between these two knots yields a new knot that also contains a valid number-pair group.

With the conclusion that adding a 1-point region between two knots containing valid number-pair groups yields a new knot that also contains a valid number-pair group, we can settle the question of the existence of Figure 4 raised in the main text. Here, however, we regard the points in Figure 4 as knots in Figure 7. The difference between Figure 7 and Figure 4 is that in Figure 4 the symbols E, F, G, H, I, J, K, L, M, N, O denote points, whereas in Figure 7 they denote knots. For Figure 7, we first regard point T together with knot M as a new knot and assign the value 1 to point T. Since the original knot M contains a valid number-pair group and all number-pairs of the new knot are a translation of the number-pairs of knot M, this new knot also contains a valid number-pair group. We then successively add 1-point regions at adjacent knots, and each addition produces a new knot containing a valid number-pair group. Finally Figure 7 is obtained. The newest knot in Figure 7 is the knot with vertex S, which also contains a valid number-pair group. Thus Figure 7 certainly exists. Returning from this viewpoint to Figure 4 in the main text, Figure 4 and Figure 7 are in fact exactly the same map, differing only in the angle from which we consider the problem; from the existence of Figure 7 we conclude that Figure 4 also exists.

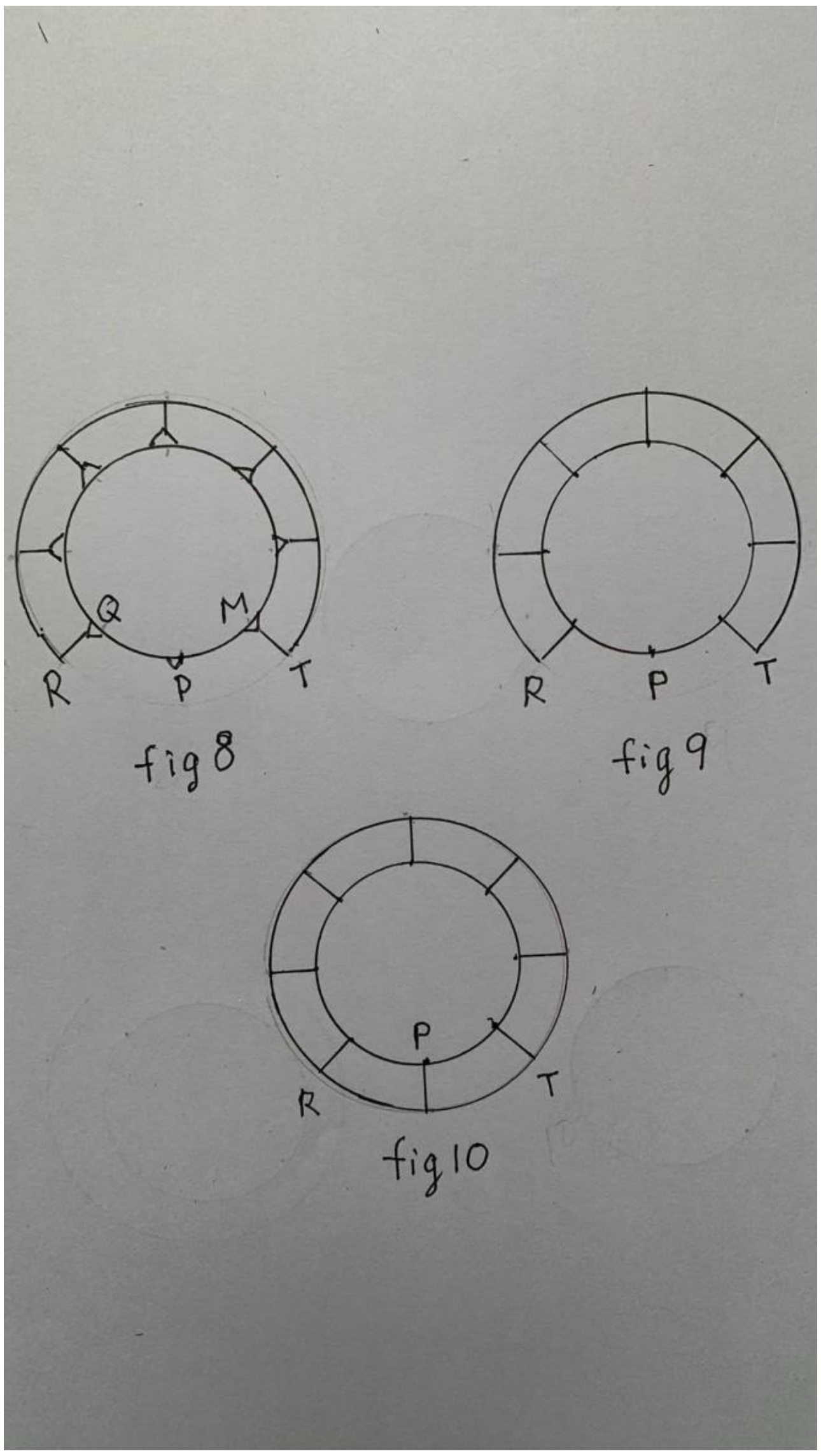

We discuss this problem further. If, on the basis of Figure 7, we continue adding 1-point regions at junction intervals, we finally reach a map with only two knots left (see Figure 8). We then change the knots in Figure 8 back into points (see Figure 9). Figures 8 and 9 are in fact exactly the same map, differing only in viewpoint. The outer boundary in Figure 9 has three intervals whose colors are pairwise different, so a 1-point region can be added only at one junction interval of Figure 9, and then a 0-point region is added, giving Figure 10. Figure 10 shows that all intervals on the original outer boundary can be colored with three colors other than some one color. If all intervals on an outer boundary can be colored with three colors, this fourth color is called the outer color of that outer boundary.

In fact, that the outer boundary of a simple closed curve has an outer color is also a necessary condition for the Four Color Theorem to hold. From the discussion above,

we obtain

**Theorem 3　An outer boundary satisfying Property A has an outer color.**

We now turn to the problem of adding simple regions to new intervals.

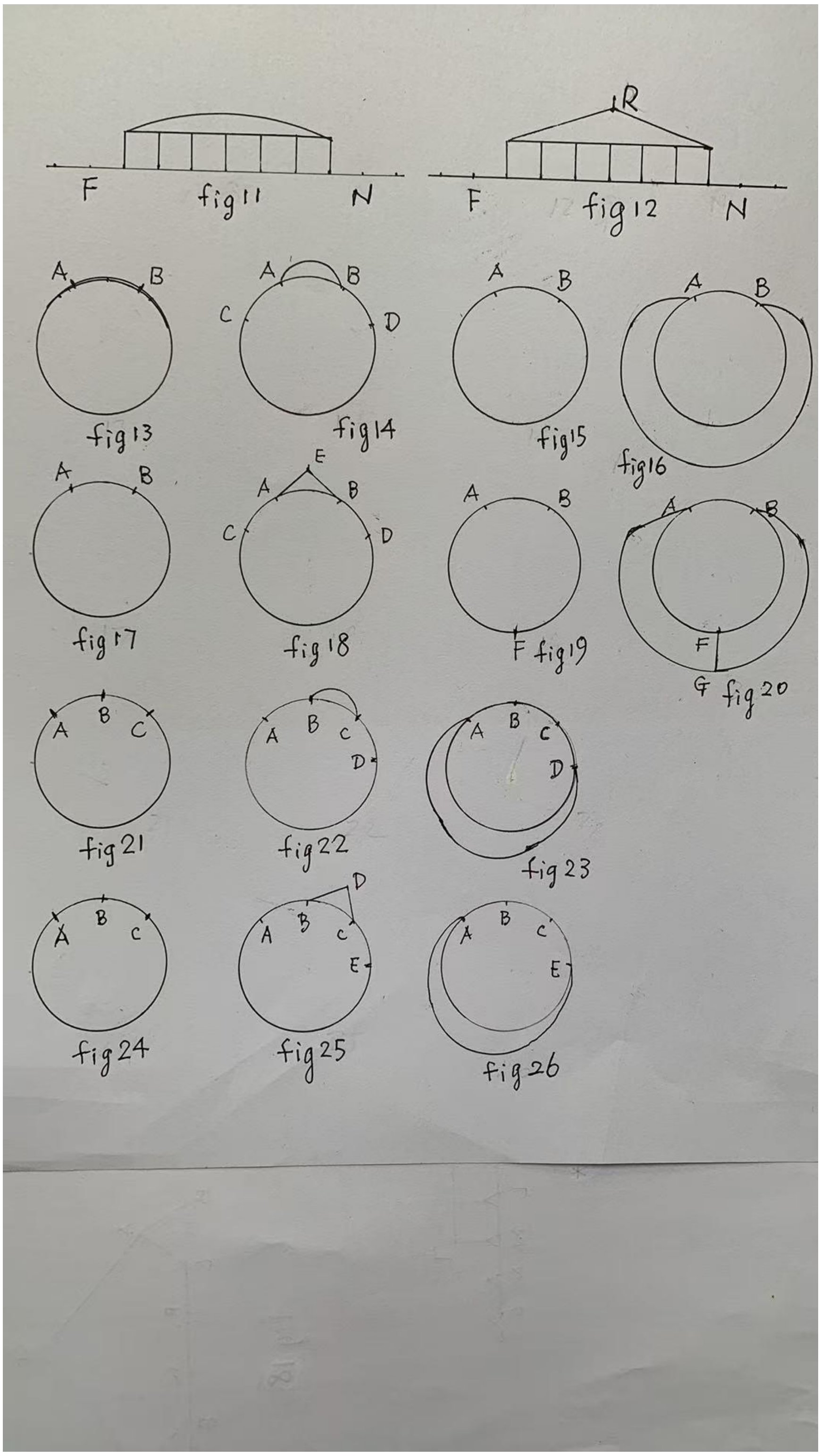

If we add a 0-point region to a new interval (see Figure 11), or add a 1-point region to a new interval (see Figure 12), we call the intervals on the outer boundary that are not old intervals boundary intervals. Figure 11 has one boundary interval FN; Figure 12 has two boundary intervals FR and RN. If a simple region is added to an old interval, since old intervals lie on an outer boundary satisfying Property A, adding the simple region does not affect whether the current outer boundary satisfies Property A; therefore we only examine the addition of simple regions to boundary intervals.

If the map has only one boundary interval and we add a 0-point region to it, the resulting new outer boundary still has only one boundary interval. If the outer boundary has one boundary interval and we add a 1-point region to it, the resulting new outer boundary has two boundary intervals. If the outer boundary has two boundary intervals and we add a 0-point region to one of them, the resulting new outer boundary has one boundary interval. If the outer boundary has two boundary intervals and we add a 1-point interval to one of them, the resulting new outer boundary has two boundary intervals. Hence the number of boundary intervals on an outer boundary is always either one or two.

We call the outer boundary containing all the old intervals the lower outer boundary, and the outer boundary containing the boundary intervals the upper outer boundary.

When the outer boundary has only one boundary interval AB (see Figure 13), we add a 0-point region to this boundary interval (see Figure 14). After adding the 0-point region we obtain a new boundary interval CD. We can now add simple regions to intervals of the lower outer boundary so that the two points A and B become the two endpoints of an interval on the lower outer boundary (see Figure 15). Then we add a 0-point region to this interval of the lower outer boundary (see Figure 16). From Figure 16 we see that, as endpoints of an interval on the upper outer boundary, the two points A and B have the same color on the boundary intervals they connect, which satisfies the earlier condition that the outer boundary has only one boundary interval.

When the outer boundary has only one boundary interval AB (see Figure 17), we add a 1-point region to this boundary interval (see Figure 18), obtaining two new boundary intervals CE and ED. We now add simple regions to intervals of the lower outer boundary so that A and B become the endpoints of two adjacent intervals on the lower outer boundary (see Figure 19). Then we add a 1-point region to interval AF of the lower outer boundary in Figure 19, and then add a 0-point region to interval GB after the 1-point region is added, giving Figure 20. From Figure 20 we see that, as endpoints of intervals on the upper outer boundary, A and B have the same color on the boundary intervals they connect, which satisfies the earlier condition that the outer boundary has only one boundary interval.

When the outer boundary has two boundary intervals AB and BC (see Figure 21), we add a 0-point region to one of them, BC, obtaining a new boundary interval AD. Then we add simple regions to intervals of the lower outer boundary so that A and D become the endpoints of an interval on the lower outer boundary (see Figure 22), and then add a 0-point region to interval AD of the lower outer boundary (see Figure 23). From Figure 23 we see that the boundary interval AB in Figure 23 and the interval BC in the figure are two adjacent intervals, so their colors are different. This satisfies the earlier condition that the outer boundary has two boundary intervals.

When the outer boundary has two boundary intervals AB and BC (see Figure 24), we add a 1-point region to one of them, BC, obtaining two new boundary intervals AD and DE (see Figure 25). Then we add simple regions to intervals of the lower outer boundary so that A and E become the endpoints of an interval on the lower outer boundary, and then add a 0-point region to interval AD of the lower outer boundary (see Figure 26). From Figure 26 we see that the interval in Figure 26 containing the interval AB of Figure 24 and the interval BC of Figure 24 are adjacent intervals, so their colors are different. This satisfies the earlier condition that the outer boundary has two boundary intervals.

We call all boundary intervals obtained from the first addition of simple regions to new intervals the first boundary-interval group; all boundary intervals obtained after adding a simple region to a boundary interval of the first group the second boundary-interval group; and in general, all boundary intervals obtained after adding a simple region to a boundary interval of the nth group the (n+1)th boundary-interval group. We continue until at last only one old interval, or one boundary interval and two old intervals, remain on the outer boundary. At that point the number of intervals on the outer boundary is less than four, so Property A is satisfied.

From the discussion above we see that if the outer boundary has the nth boundary-interval group, we can always add simple regions on the lower outer boundary so that afterward the requirement for the (n-1)th boundary-interval group is satisfied. Hence we can always add simple regions step by step on old intervals, satisfying the conditions of all boundary-interval groups until the requirement of the first boundary-interval group is met.

Here every boundary-interval group is obtained by continuously adding simple regions to old intervals in Figure 11 or Figure 12, so the outer boundary at that point satisfies Property A. In particular, the outer boundary containing the old intervals from the initial outer boundary when the first boundary group is added in the first boundary-interval group also satisfies Property A. By Theorem 3, the outer boundary at that point has an outer color. Therefore the map before adding simple regions in Figure 11 or Figure 12 exists, and we know that the outer boundaries of both Figure 11 and Figure 12 satisfy Property A.

**Theorem 4** If an n-point region is added to any interval of an outer boundary satisfying Property A, the resulting new outer boundary also satisfies Property A. (n-Point Region Theorem)

**Proof** If the outer boundary has only two or three intervals and we add an n-point region to any interval, then since no old intervals remain after the n-point region is added in this case, it can also be regarded as adding a 0-point or 1-point region to an interval of a polygonal outer boundary, the color of that polygon being the color of the added n-point region. Hence the result satisfies Property A.

When the outer boundary has more than three intervals, we already know that adding a 0-point, 1-point, or 2-point region to any interval of an outer boundary satisfying Property A yields a new outer boundary satisfying Property A. Now assume that adding a region with fewer than m points to any interval of an outer boundary satisfying Property A yields a new outer boundary satisfying Property A. Consider the case of adding an m-point region.

If a simple region is added to a new interval of the new outer boundary, the number of new points decreases, which amounts to adding a region with fewer than m points to an outer boundary satisfying Property A; by the induction hypothesis, the new outer boundary satisfies Property A.

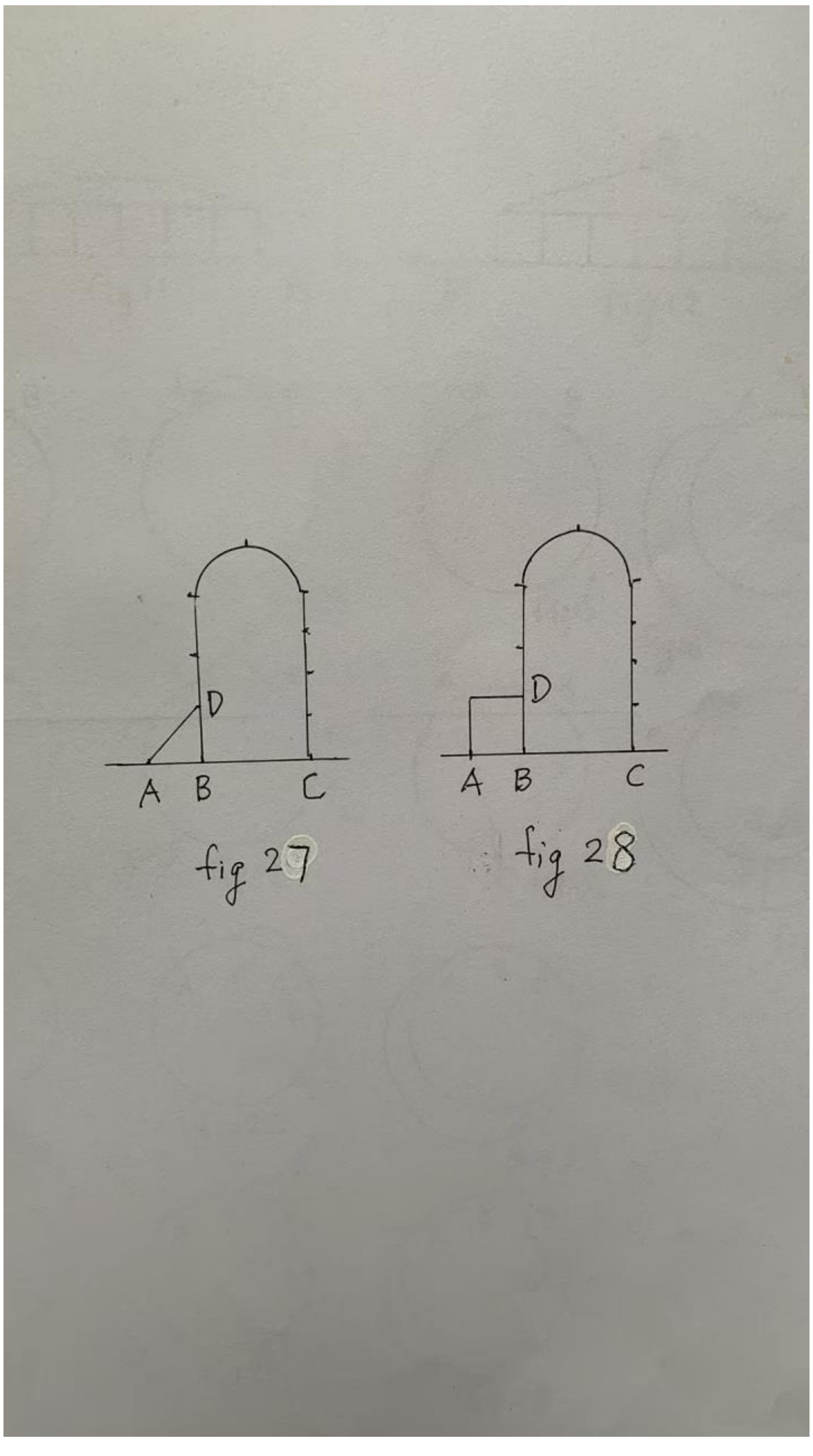

If a 0-point region is added to the junction interval AD of the new outer boundary (see Figure 27), Figure 27 is obtained by first adding an m-point region to interval BC of the original outer boundary and then adding a 0-point region to its junction interval AD. It can also be regarded as obtained by first adding a 1-point region to interval AB of the original outer boundary and then adding an (m-1)-point region to junction interval DC of the new outer boundary. Because the original outer boundary satisfies Property A, the outer boundary after the 1-point region is added still satisfies Property A; then adding an (m-1)-point region to that outer boundary, by the induction hypothesis, yields a new outer boundary that also satisfies Property A.

If a 1-point region is added to the junction interval AD of the new outer boundary (see Figure 28), Figure 28 is obtained by first adding an m-point region to interval BC of the original outer boundary and then adding a 1-point region to junction interval AD of the new outer boundary. It can also be regarded as obtained by first adding a 2-point region to interval AB of the original outer boundary and then adding an (m-1)-point region to junction interval DC of the new outer boundary. By Theorem 2, the outer boundary after adding a 2-point region to the original outer boundary satisfies Property A; then adding an (m-1)-point region to junction interval DC of that outer boundary, by the induction hypothesis, gives the newest outer boundary satisfying Property A.

If simple regions are added only to old intervals, since old intervals lie on the initial outer boundary satisfying Property A, the problem remains the old one: does adding an m-point region to some interval of an outer boundary satisfying Property A yield a new outer boundary that still satisfies Property A? Only the number of old intervals has decreased. Since the number of old intervals is finite, it will eventually reach zero. With no old intervals, the colors of the intervals on the new outer boundary correspond one-to-one to the colors of the intervals on a polygonal outer boundary, the color of that polygon being the color of the added n-point region; hence the outer boundary at that point satisfies Property A.

In summary, after adding an m-point region to an outer boundary satisfying Property A, no matter how simple regions are subsequently added, one can always explicitly obtain an outer boundary satisfying Property A.

Thus by induction, after adding an arbitrary n-point region to an outer boundary satisfying Property A, its outer boundary still satisfies Property A.

That is, Theorem 4 holds.

### Coloring the Island with No Interior Regions

Before coloring, we first state a fact: for any simple closed curve in a map with no interior regions, taken as an outer boundary, one can always find an interval on it such that adding to this interval a region outside the boundary that has this interval as part of it is precisely what we previously called adding an n-point region to the interval, that is, apart from this interval, no other point of the region lies on the simple closed curve.

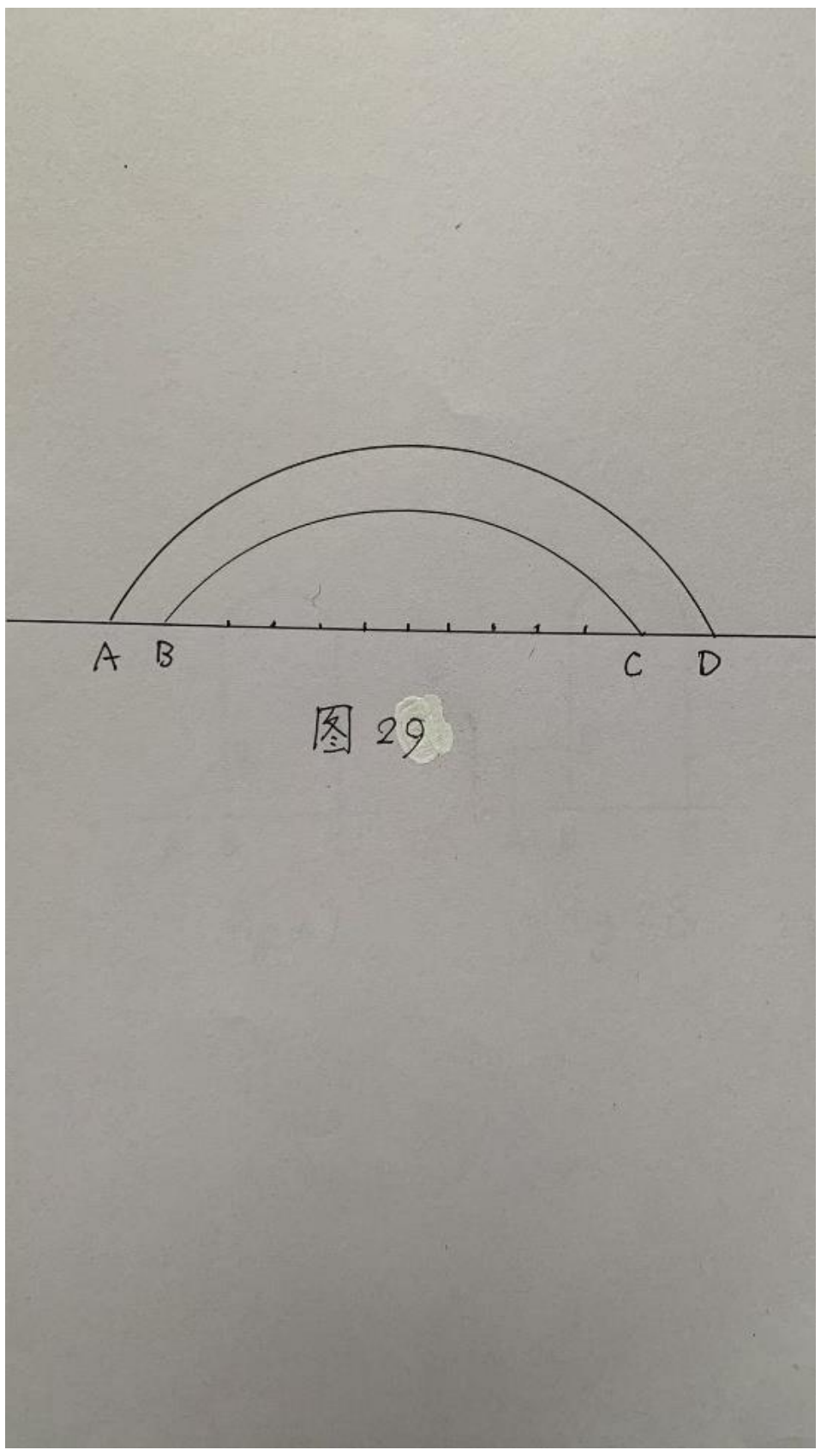

See Figure 29. If the region on interval AB of some outer boundary in the given map has another edge CD lying on that outer boundary, then clearly AB and CD are not adjacent and there are further intervals between them. If among these intervals we can find two that are two edges of the same region, we may continue our examination within these two intervals, until the regions on all intervals between some two intervals have no points other than the interval itself lying on the outer boundary. Since the number of intervals between two intervals is finite, this result can always be achieved: the regions added on these intervals are all what we called adding an n-point region to the interval.

We now proceed to color the simple island.

Since the sea is already colored, it is a colored polygon with respect to this island. By Theorem 1, the outer boundary of the sea satisfies Property A. By the fact just stated, we can always find an interval on this outer boundary and add to it a region outside the boundary, a region that has no part other than the chosen interval lying on the original outer boundary, that is, we add an n-point region. By Theorems 1, 2, 3, and 4, the outer boundary after this region is added still satisfies Property A. This operation can be continued until all regions of the simple island have been added. Note that an added region is already colored at the same time, so after all regions of the island have been added, we have in fact colored all regions of the island with no interior regions using four colors.

We stated at the beginning of this paper that if any simple island in the sea can be colored, then any given map can be colored, that is, the Four Color Theorem holds.

**Discussion and Conclusion**

The method proposed in this paper provides a constructive, non-enumerative framework for proving the Four Color Theorem. Its core advantages are as follows:
Clarity of concepts: By choosing an arbitrary region in any given spherical map as the sea, then defining simple islands of the sea, we show that as long as every simple island in the sea can be colored with four colors, any given map is guaranteed to be colorable. By defining outer boundaries, the region-adding operation, and Property A on spherical maps, we decompose the complex global coloring problem into a series of manageable local steps.
Rigor of logic: The four theorems form a complete logical chain, ensuring the sustainability and correctness of the coloring process on spherical maps.